\documentclass[12pt,a4]{article}
\usepackage{amssymb,amsmath,graphicx,enumerate}
\usepackage{longtable,mathrsfs}
\usepackage{hyperref}
\graphicspath {P:\Dole}

\DeclareMathAlphabet{\mathpzc}{OT1}{pzc}{m}{it}

\numberwithin{equation}{section}

\newtheorem{theorem}{Theorem}[section]
\newtheorem{lemma}{Lemma}[section]

\newtheorem{definition}{Definition}[section]
\newtheorem{example}{Example}[section]

\title{Linear Canonical Transform On Boehmian Space}
\author{S. K. Panchal and Pravinkumar V. Dole\\
 Department of Mathematics, \\ Dr. Babasaheb Ambedkar Marathwada University,\\ Aurangabad-431004 (M.S.) India.\\
E-mail ID- drpanchalsk@gmail.com, pvasudeo.dole@gmail.com}
\begin{document}
\maketitle{\bf Abstract:}
The aim of this paper is to constructs Boehmian space, the linear canonical transform for Boehmians is define and to study its properties.


\section{Introduction}
\indent\qquad The most recent generalizations of functions is the theory of Boehmians. The idea of construction of Boehmians was initiated by the concept of regular operators introduced by Boehme \cite{B}. Regular operators form a subalgebra of the field of
Mikusinski operators and they include only such functions whose support is bounded from the left. In a concrete case, the space of Boehmians contains all regular operators, all distributions and some objects which are neither operators nor distributions. The space of Boehmians is the new class of generalized functions which opened the new door to area of research in mathematics. The construction of Boehmians is given by Mikusinski and Mikusinski \cite{JP, PM1, PM3, PM4}. Mikusinski and Nemzer has studied Fourier and Laplace transform for Boehmians in \cite{PM2} and \cite{N} respectively. Zayed \cite{AZ} extended the fractional Fourier transform to class of integrable Boehmians. Singh studied fractional integrals of fractional Fourier transform for integrable Boehmians in \cite{SB}. The Fourier, Laplace and fractional Fourier transforms are the special cases of linear canonical transform (LCT) and has many applications in several areas, like signal processing and optics \cite{BRY}. This lead to study linear canonical transform for integrable Boehmians in \cite {PS1}. In one hand we constructs Boehmian spaces and other hand defined linear canonical transform for Boehmians. Further, we obtain its properties like one-to-one, onto, continuous from one Boehmian space to another Boehmian space and other basic properties in the space of Boehmians.

The linear canonical transform of real valued function $f$ is defined \cite{DQY, MQ} as,
\begin{align}
\mathcal{L}_{A}[f(t)](u)=F_{A}(u)=\left\{ {\begin{array}{*{20}{l}}{\sqrt{\frac{1}{2\pi ib}}\int_{-\infty}^{\infty}e^{\frac{i}{2}[\frac{a}{b}t^{2}-\frac{2}{b}ut+(\frac{d}{b})u^{2}]}f(t)dt\,\, for\,\, b\neq 0,}\\
{\sqrt{d}e^{\frac{i}{2}cdu^{2}} f(du) \qquad\qquad\qquad\qquad for\,\, b=0,}
\end{array}} \right.
\end{align}
where $\mathcal{L}_{A}$ is the unitary linear canonical transform operator with parameter $A=(a, b, c, d)$, $a,b,c,d$ are real number satisfying $ad-bc=1$. The inverse transform for linear canonical transform is given by a linear canonical transform having the parameter $A^{-1}=(d, -b, -c, a)$ and $\mathcal{L}_{A^{-1}}$ is the inverse LCT operator. For value of parameters as, $a=cos\theta, b=sin\theta, c=-sin\theta, d=cos\theta $ then LCT become fractional Fourier transform, in particular, when $\theta= \frac{\pi}{2}$ then LCT become Fourier transform and $a=0, b=i, c=i, d=0 $ then LCT becomes Laplace transform.

 Let $\mathcal{L}^{1}(\mathbb{R})$ be the space of all complex valued absolutely integrable functions on $\mathbb{R}$ with norm $||f||_{1}=\int_{\mathbb{R}}|f(t)|dt \leq M_{1}$ and $\mathcal{L}^{2}(\mathbb{R})$ be the space of all complex valued absolutely square integrable functions on $\mathbb{R}$ with norm $||g||_{2}=\big{(}\int_{\mathbb{R}}|g(t)|^{2}dt\big{)}^{\frac{1}{2}} \leq M_{2}$, for some $M_{1}, M_{2} > 0$. Let $\mathcal{L}^{1}(\mathbb{R})\cap \mathcal{L}^{2}(\mathbb{R})$ is denoted by $\mathcal{L}^{1, 2}(\mathbb{R})$.
\begin{definition}\cite{Pathak}(Regular Distributions) Let $f$ be the locally integrable function, i.e. absolutely integrable on every finite interval on $\mathbb{R}$, then distribution generated by $f$ is called  regular distributions.
\end{definition}
 We see that $f\in\mathcal{L}^{1, 2}(\mathbb{R})$ then $\mathcal{L}_{A}(f)$ and $\mathcal{L}_{A^{-1}}(f)$ are the members of $\mathcal{L}^{1, 2}(\mathbb{R})$.
\begin{definition} \cite{DQY} Let the weight function $W(t, \tau)= e^{i\tau(\tau-t)\frac{a}{b}}$. For any two function $f$ and $g$ the convolution operation $*^{A}$ is defined as,
\begin{align}
h(t)=(f *^{A} g)(t)=\int_{-\infty}^{\infty}f(\tau)g(t-\tau)W(t,\tau)d\tau
\end{align}
\end{definition}
\begin{theorem} \cite{DQY} (New Convolution Theorem)\label{thm1.1}\\
Let $h(t)=(f *^{A} g)(t)$ and $ H_{A}(u), F_{A}(u), G_{A}(u)$ denote the linear canonical transform of $h(t), f(t)$ and $g(t)$ respectively, then
\begin{align}
H_{A}(u)=\sqrt{2i\pi b}\, e^{-i(\frac{d u^{2}}{2b})}F_{A}(u)G_{A}(u).
\end{align}
\end{theorem}
\section{Preliminary Results}
\indent\qquad In this section we obtain some results which are require to construct the Boehmian space. 
\begin{lemma}\label{lemma2.1}
 Let $f\in \mathcal{L}^{1}(\mathbb{R})$ and $g\in \mathcal{L}^{2}(\mathbb{R})$ then the $(f*^{A} g)$ is in $\mathcal{L}^{2}(\mathbb{R})$.
\end{lemma}
\begin{lemma}
 The space $(\mathcal{L}^{1, 2}(\mathbb{R}), *^{A})$ is commutative semi group.
\end{lemma}
\begin{theorem}(Plancherel type theorem) Let the sequence $f_{n}\in \mathcal{L}^{1, 2}(\mathbb{R})$ and $f_{n}\rightarrow f$ on $\mathcal{L}^{2}(\mathbb{R})$ then  $\mathcal{L}_{A}(f_{n})\rightarrow \mathcal{L}_{A}(f)$ in $\mathcal{L}^{2}(\mathbb{R})$ as $n\rightarrow \infty$.
\end{theorem}
\begin{definition}
Analogous to Plancherel type theorem for $f\in \mathcal{L}^{2}(\mathbb{R})$, we define $\mathcal{L}_{A}(f)$ by $\mathcal{L}^{2}-\lim_{n\rightarrow \infty} \mathcal{L}_{A}(f_{n})$, where $f_{n}\in \mathcal{L}^{1, 2}(\mathbb{R})$.
\end{definition}
Let $\bigtriangledown$ be the set of all sequences of continuous real functions $\{\delta_{n}\}$ from $\mathcal{L}^{1, 2}(\mathbb{R})$ having compact support on $\mathbb{R}$ with the following properties:
\begin{enumerate}
\item[(i)] $\quad\int_{\mathbb{R}}e^{i\frac{at^{2}}{2b}}\delta_{n}(t)dt=1$,  $\forall \, n\in \mathbb{N}$,
\item[(ii)] $\quad\lim_{n\rightarrow \infty}\int_{|t|>\epsilon} |\delta_{n}(t)| dt = 0$ for each $\epsilon >0$.\end{enumerate}
 The members of $\bigtriangledown$ are called {\it delta sequences}.
\begin{example}
Let $a, b\in\mathbb{R}; b\neq 0$, consider the sequence
\begin{align*}
\delta_{n}(t)=\left\{ {\begin{array}{*{20}{l}}{e^{-i\frac{at^{2}}{2b}}t \qquad\qquad\qquad for \quad 0\leq t\leq \frac{1}{n},}\\{e^{-i\frac{at^{2}}{2b}}n^{2} (\frac{2}{n}-t)\qquad\quad\,\, for \quad \frac{1}{n}\leq t\leq \frac{2}{n},}\\
{0 \qquad\qquad\qquad\qquad\qquad\quad otherwise}.
\end{array}} \right.
\end{align*}
\end{example}
\begin{lemma}\label{lemma2.3}
  Let $\{\varphi_{n}\}, \{\psi_{n}\}\in \bigtriangledown$ then $(\varphi_{n}*^{A}\psi_{n}) \in \bigtriangledown$ for all $n\in \mathbb{N}$.
\end{lemma}
\begin{lemma}
Let $f\in \mathcal{L}^{1, 2}(\mathbb{R})$ and $\{\psi_{n}\}\in \bigtriangledown $ then $f *^{A}\psi_{n}\rightarrow f$ as $n\rightarrow \infty$ in $\mathcal{L}^{2}(\mathbb{R})$.
\end{lemma}
\section{LCT For Boehmians}
 \indent\qquad A pair of sequences $(f_{n}, \varphi_{n})$ is called a quotient of the sequences, denoted by $f_{n}/\varphi_{n}$, where each $n\in\mathbb{N}$, $f_{n}\in \mathcal{L}^{1, 2}(\mathbb{R})$ and $\{\varphi_{n}\}\in \bigtriangledown $ such that $f_{m}*^{A}\varphi_{n}=f_{n}*^{A}\varphi_{m}$ holds $\forall\, m,n\in \mathbb{N}$. Two quotients of sequences $f_{n}/ \varphi_{n}$ and $g_{n} / \psi_{n}$ are equivalent if $f_{n}*^{A}\psi_{n}=g_{n}*^{A}\varphi_{n}$ $\forall\, n\in\mathbb{N}$. This is an equivalence relation. The equivalence class of quotient of sequences is called a {\it Boehmian}. The space of all Boehmians is denoted by $\mathcal{B}_{\mathcal{L}^{1, 2}}=\mathcal{B}_{\mathcal{L}^{1, 2}}(\mathcal{L}^{1, 2}(\mathbb{R}), \bigtriangledown, *^{A})$ and the members of $\mathcal{B}_{\mathcal{L}^{1, 2}}$ are denoted by $F=[f_{n}/\varphi_{n}]$. The function $f\in \mathcal{L}^{1, 2}(\mathbb{R})$ can be identified with the Boehmian $[(f *^{A}\delta_{n})/\delta_{n}]$, where $\{\delta_{n}\}$ is the delta sequence. Let $F=[f_{n}/\varphi_{n}]$, then $F*^{A}\delta_{n}=f_{n}\in \mathcal{L}^{1, 2}(\mathbb{R})$ $\forall\, n\in\mathbb{N}$.
\begin{definition}
  A sequence of Boehmians $F_{n}$ is called $\Delta-$convergent to a Boehmian $F$ ($\Delta-\lim F_{n}=F$) if there exist a delta sequence $\{\delta_{n}\}$ such that  $(F_{n}-F)*^{A}\delta_{n}\in \mathcal{L}^{1, 2}(\mathbb{R})$, for every $n\in\mathbb{N}$ and that $\|(F_{n}-F)*^{A}\delta_{n}\|_{2}\rightarrow 0$ as $n\rightarrow \infty$.
\end{definition}
\begin{definition}
A sequence of Boehmians $F_{n}$ is called $\delta-$convergent to a Boehmian $F$ ($\delta-\lim F_{n}=F$) if there exist a delta sequence $\{\delta_{n}\}$ such that  $F_{n}*^{A}\delta_{k}\in \mathcal{L}^{1, 2}(\mathbb{R})$ and $F *^{A} \delta_{k}\in \mathcal{L}^{1, 2}(\mathbb{R})$ for every $n, k\in\mathbb{N}$ and that $\|(F_{n}-F)*^{A}\delta_{k}\|_{2}\rightarrow 0$ as $n\rightarrow \infty$ for each $k\in\mathbb{N}$. 
\end{definition}
\par Let $\{\delta_{n}\}$ is a delta sequence, then $\delta_{n}/\delta_{n}$ represents an Boehmian. Since the Boehmian $[\delta_{n}/\delta_{n}]$ corresponds to Dirac delta distribution $\delta$, all the derivative of $\delta$ are also Boehmian. If $\{\delta_{n}\}$ is infinitely differentiable and bounded support, then the $k^{th}$ derivative of $\delta$ is define by $\delta^{(k)}=[\delta_{n}^{(k)}/\delta_{n}]\in \mathcal{B}_{\mathcal{L}^{1, 2}}$, for each $k\in\mathbb{N}$. The $k^{th}$ derivative of Boehmian $F\in \mathcal{B}_{\mathcal{L}^{1, 2}}$ is define by $F^{(k)}=F *^{A} \delta^{(k)}$.The scalar multiplication, addition and convolution in $\mathcal{B}_{\mathcal{L}^{1, 2}}$ are define as,
\begin{align*}
\lambda[f_{n}/\varphi_{n}]&=[\lambda f_{n}/\varphi_{n}]\\
[f_{n}/\varphi_{n}]+[g_{n}/\psi_{n}]&=[(f_{n}*^{A}\psi_{n}+g_{n}*^{A}\varphi_{n})/\varphi_{n}*^{A}\psi_{n}]\\
[f_{n}/\varphi_{n}]*^{A}[g_{n}/\psi_{n}]&=[(f_{n}*^{A} g_{n})/(\varphi_{n}*^{A}\psi_{n})].
\end{align*}
\begin{lemma} Let $\Delta-\lim F_{n}= F$ in $\mathcal{B}_{\mathcal{L}^{1, 2}}$ then $\Delta-\lim F_{n}^{(k)}= F^{(k)}$ for $\forall  \, k\in\mathbb{N}$ in $\mathcal{B}_{\mathcal{L}^{1, 2}}$.
\end{lemma}

Let $\bigtriangledown_{0}=\{\mathcal{L}_{A}(\delta_{n}); \{\delta_{n}\}\in \bigtriangledown \}$ be the space of complex valued functions on $\mathbb{R}$, the operation $\cdot$ is pointwise multiplication and $C_{0}(\mathbb{R})$ be the space of all continuous functions vanishing at infinity on $\mathbb{R}$ then we construct the another space of Boehmians, denoted by $\mathcal{B}_{\bigtriangledown}= \mathcal{B}_{\bigtriangledown}(\mathcal{L}^{2}(\mathbb{R}), C_{0}(\mathbb{R})\cap \mathcal{L}^{2}(\mathbb{R}), \cdot, \bigtriangledown_{0})$. This is the range of linear canonical transform on $\mathcal{B}_{\mathcal{L}^{1, 2}}$ and each element of $\mathcal{B}_{\bigtriangledown}$ is denoted by $\mathcal{L}_{A} (f_{n})/ \mathcal{L}_{A}(\delta_{n})$ for all $n\in\mathbb{N}$, where $\{f_{n}\}\in \mathcal{L}^{1, 2}(\mathbb{R}) $.
\begin{lemma} Let $f, g\in \mathcal{L}^{2}(\mathbb{R}); \varphi, \psi \in C_{0}(\mathbb{R})$ and $\lambda\in\mathbb{C}$ then\\ (i) $f\cdot\varphi \in \mathcal{L}^{2}(\mathbb{R})$\\ (ii) $(f+g)\cdot\varphi =f\cdot\varphi + f\cdot\varphi$\\ (iii) $(\lambda f)\cdot\varphi=\alpha(f\cdot\varphi)$\\ (iv) $f\cdot(\varphi\cdot\psi)=(f\cdot\varphi)\cdot\psi$.
\end{lemma}
\begin{lemma} Let $f_{n}\rightarrow f$ as $n\rightarrow\infty$ in $\mathcal{L}^{2}(\mathbb{R})$ and $\varphi\in C_{0}(\mathbb{R})$ then $f_{n}\cdot\varphi\rightarrow f\cdot\varphi$ in $\mathcal{L}^{2}(\mathbb{R})$.
\end{lemma}
\begin{lemma}\label{lemma3.7}
Let $\{\delta_{n}\}\in \bigtriangledown$ then $\mathcal{L}_{A}(\delta_{n})$ converges uniformly on each compact set to a constant function $1$ in $\mathcal{L}^{2}(\mathbb{R})$.
\end{lemma}
\begin{lemma} Let $f_{n}\longrightarrow f$ as $n\longrightarrow\infty$ in $\mathcal{L}^{1, 2}(\mathbb{R})$ and $\mathcal{L}_{A}(\varphi_{n})\in \bigtriangledown_{0}$ then $f_{n}\cdot \mathcal{L}_{A}(\varphi_{n})\rightarrow f$ in $\mathcal{L}^{2}(\mathbb{R})$.
\end{lemma}
\begin{lemma} Let $\mathcal{L}_{A}(\varphi_{n}), \mathcal{L}_{A}(\psi_{n})\in \bigtriangledown_{0}$ then $\mathcal{L}_{A}(\varphi_{n})\cdot \mathcal{L}_{A}(\psi_{n}) \in \bigtriangledown_{0}$.
\end{lemma}
\textbf{Proof:} Let $ \mathcal{L}_{A}(\varphi_{n}), \mathcal{L}_{A}(\psi_{n}) \in C_{0}(\mathbb{R})$ From theorem \eqref{thm1.1} and lemma \eqref{lemma2.3} we get
$\mathcal{L}_{A}(\varphi_{n})\cdot \mathcal{L}_{A}(\psi_{n})=\frac{e^{\frac{i}{2}(\frac{d}{b})u^{2}}}{\sqrt{2\pi ib}}\mathcal{L}_{A}(\varphi_{n}*^{A}\psi_{n})\in\bigtriangledown_{0}$.$\hfill\blacksquare$
\begin{definition}
Let $\{f_{n}\}\in \mathcal{L}^{1, 2}(\mathbb{R})$ and $\{\delta_{n}\}\in\bigtriangledown$, we define the linear canonical transform $\mathcal{L}_{A} : \mathcal{B}_{\mathcal{L}^{1, 2}}\longrightarrow\mathcal{B}_{\bigtriangledown}$ as
 \begin{align}
 \mathcal{L}_{A} [f_{n}/\delta_{n}]= \mathcal{L}_{A} (f_{n})/ \mathcal{L}_{A} (\delta_{n}) \qquad for \quad [f_{n}/\delta_{n}]\in \mathcal{B}_{\mathcal{L}^{1, 2}}.
\end{align}
\end{definition}
\indent\quad The linear canonical transform  on $\mathcal{B}_{\mathcal{L}^{1, 2}}$ is well defined. Indeed if  $[f_{n}/\delta_{n}]\in \mathcal{B}_{\mathcal{L}^{1, 2}}$, then $f_{n}*^{A}\delta_{m}=f_{m}*^{A}\delta_{n}$ for all $m, n \in \mathbb{N}$. Applying the linear canonical transform on both sides, we get $\mathcal{L}_{A} (f_{n}) \mathcal{L}_{A} (\delta_{m})=\mathcal{L}_{A} (f_{m}) \mathcal{L}_{A} (\delta_{n})$ for all $m, n \in \mathbb{N}$ and hence $\mathcal{L}_{A} (f_{n})/ \mathcal{L}_{A} (\delta_{n})\in \mathcal{B}_{\bigtriangledown}$. Further if $[f_{n}/\psi_{n}]=[g_{n}/\delta_{n}]\in \mathcal{B}_{\mathcal{L}^{1, 2}}$ then we have $f_{n}*^{A}\delta_{n}=g_{n}*^{A}\psi_{n}$ for all $ n \in \mathbb{N}$. Again applying the linear canonical transform on both sides, we get $\mathcal{L}_{A} (f_{n}) \mathcal{L}_{A} (\delta_{n})=\mathcal{L}_{A} (g_{n}) \mathcal{L}_{A} (\psi_{n})$ for all $ n \in \mathbb{N}$. i.e. $\mathcal{L}_{A} (f_{n})/ \mathcal{L}_{A} (\psi_{n})=\mathcal{L}_{A} (g_{n})/ \mathcal{L}_{A} (\delta_{n})$ in $\mathcal{B}_{\bigtriangledown}$.
\begin{lemma}\label{2.3.7}
Let $[f_{n}/ \varphi_{n}]\in B_{\mathcal{L}^{1, 2}}$ then the linear canonical transform of the sequence
\begin{align}
\mathcal{L}_{A} [f_{n}](u)= \sqrt{\frac{1}{2\pi ib}}e^{\frac{i}{2}(\frac{d}{b})u^{2}}\int_{-\infty}^{\infty}e^{\frac{-i}{b}ut}e^{\frac{i}{2}\frac{a}{b}t^{2}}f_{n}(t)dt
\end{align}
converges uniformly on each compact set in $\mathbb{R}$.
\end{lemma}
\begin{definition} In view of the above proof of lemma \eqref{2.3.7}, the linear canonical transform of Boehmian in the space of continuous functions on $\mathbb{R}$ is define as,
\begin{align*}
\mathcal{L}_{A}[F]=\lim_{n\rightarrow\infty}\mathcal{L}_{A} (f_{n}).
\end{align*}
\end{definition}
\begin{theorem}
The linear canonical transform $\mathcal{L}_{A}:\mathcal{B}_{\mathcal{L}^{1, 2}}\longrightarrow \mathcal{B}_{\bigtriangledown}$ is consistent with $\mathcal{L}_{A}: \mathcal{L}^{2}(\mathbb{R})\longrightarrow \mathcal{L}^{2}(\mathbb{R})$.
\end{theorem}
\begin{theorem}
The linear canonical transform $\mathcal{L}_{A}:\mathcal{B}_{\mathcal{L}^{1, 2}}\longrightarrow \mathcal{B}_{\bigtriangledown}$ is a bijection.
\end{theorem}
\begin{theorem}
Let $F, G \in \mathcal{B}_{\mathcal{L}^{1, 2}}$ then
\begin{enumerate}
\item[(a)] $\quad\mathcal{L}_{A}[F + \lambda G]= \mathcal{L}_{A} (F)+ \lambda  \mathcal{L}_{A} (G)$, for any complex $\lambda$.
\item[(b)] $\quad\mathcal{L}_{A}[e^{ikt} F](u)=e^{\frac{-idk(2u-bk)}{2}}\mathcal{L}_{A}[F](u-bk)$, for $k\in\mathbb{R}$.
\item[(c)] $\quad\mathcal{L}_{A} [F (t+\tau)](u)=e^{i(2u+a\tau)\frac{\tau}{2b}}\mathcal{L}_{A}[e^{\frac{-ia}{b}x\tau}F(x)](u)$.
\item[(d)] $\quad\mathcal{L}_{A} [F^{(2)}](u)=\bigg{[}\bigg{(}\frac{iu}{b}\bigg{)}^{2}+\frac{ia}{b}\bigg{]}\mathcal{L}_{A} [F(t)](u).$
\end{enumerate}
\end{theorem}
\begin{theorem}
 Let $F, G\in \mathcal{B}_{\mathcal{L}^{1, 2}}$ then $\mathcal{L}_{A} (F *^{A} G)= \mathcal{L}_{A}(F) \mathcal{L}_{A}(G)$.
\end{theorem}
\begin{theorem}
Let $\delta-\lim F_{n}=F$ for $F_{n}, F\in \mathcal{B}_{\mathcal{L}^{1, 2}}$ then $\mathcal{L}_{A} (F_{n})\rightarrow \mathcal{L}_{A} (F)$ uniformly on each compact set of $\mathbb{R}$.
\end{theorem}
\textbf{Proof:} Let $\{\delta_{m}\}$ be a delta sequence such that $F_{n}*^{A}\delta_{m}, F*^{A}\delta_{m}\in \mathcal{L}^{1, 2}(\mathbb{R})$ for all $n,m\in \mathbb{N}$ and $\Vert(F_{n}-F)*^{A}\delta_{m}\Vert_{2}\rightarrow 0$ as $n\rightarrow \infty$ for each $m\in\mathbb{N}$. Let $M$ be a compact set in $\mathbb{R}$ then $\mathcal{L}_{A}(\delta_{m})>0$ on $M$ for all most $m\in\mathbb{N}$. Since $\mathcal{L}_{A}(\delta_{m})$ is a continuous function and $\mathcal{L}_{A}(F_{n})*^{A} \mathcal{L}_{A}(\delta_{m})- \mathcal{L}_{A}(F)*^{A} \mathcal{L}_{A}(\delta_{m})=((\mathcal{L}_{A}(F_{n})- \mathcal{L}_{A}(F))*^{A} \mathcal{L}_{A}(\delta_{m}))$, implies $\Vert (\mathcal{L}_{A}(F_{n})- \mathcal{L}_{A}(F))*^{A} \mathcal{L}_{A}(\delta_{m}) \Vert_{2}\rightarrow 0$ as $n\rightarrow \infty$ for each $m\in\mathbb{N}$. Thus $\mathcal{L}_{A}(F_{n})\rightarrow \mathcal{L}_{A}(F)$ uniformly on $M$. $\hfill\blacksquare$


\begin{thebibliography}{99}

\bibitem{B} T. K. Boehme; The support of Mikusinski operators, {\it Trans. Amer. Math. Soc.}, 176, 319-334, (1973).

\bibitem{BRY} Deng Bing, Tao Ran and Wang Yue. Convolution theorems for the linear canonical transform and their applications. {\it Sci. China Series F: Inf. Sci.}, 49(5), 592-603, (2006).

\bibitem{DQY} Deyun Wei, Qiwen Ran and Yong Li; New convolution theorem for the linear canonical transform and its translation invariance property, {\it Optik} 123, 1478-1481, (2012).

\bibitem{PS1} Pravinkumar V. Dole and S. K. Panchal, Linear canonical transform for Integrable Boehmians, {\it Int. J. Pure Appl. Math.}, 116(1), 91-96, (2017).

\bibitem{JP} J.Mikusinski and P.Mikusinski; Quotients de suites et leurs applications dans l’analyse
fonctionnelle, {\it C.R. Acad. Sci. Paris Ser. I Math.}, 293, 463-464, (1981).

\bibitem{PM1} P. Mikusinski; Convergence of Boehmianes, {\it Japan. J. Math}, 9(1), 169-179,(1983).

\bibitem{PM2} P. Mikusinski; Fourier transform for integrable Boehmians, {\it Rocky Moun-tain J. Math.}, 17(3), 577-582,(1987).

\bibitem{PM3} P. Mikusinski; Boehmians and generalized functions, {\it Acta. Math. Hungarica}, 51, 159-179, (1988).

\bibitem{PM4} P. Mikusinski; Transform of Boehmians, {\it Different Aspects of Differentiability Dissertationes Mathematicae},  340, 201-206,(1995).

\bibitem{MQ} M. Moshinsky and C. Quesne; Linear canonical transformations and their unitary representations, {\it J. Math. Phys}, 12(8), 1772-1783,(1971).

\bibitem{N} Dennis Nemzer; Laplace transforms on a class of Boehmians, {\it Bull. Austral. Math. Soc.}, 46, 347-352, (1992).

\bibitem{Pathak} R. S. Pathak, \emph{A Course in Distributional Theory and Applications}, Narosa Publication House, New Delhi 2001.

\bibitem{SB} A. Singh and P. K. Banergi; Fractional integrals of fractional Fourier transform for integrable Boehmians,
{\it Proceedings of the National Academy of Sciences, India Section A: Physical Sciences}, 2017. 

\bibitem{W} Walter Rudin; {\it Real and Complex Analysis}, Third Edition, McGraw-Hill, New York, 1987.

\bibitem{AZ} A. I. Zayed; Fractional Fourier transform of generalized functions, {\it Integ. Trans. Spl. Funct.},7, 299-312, (1998).

\end{thebibliography}
\end{document}